\documentclass[final]{siamltex}
\usepackage{graphicx,bm,amssymb,amsmath,hyperref}

\newcommand{\bi}{\begin{itemize}}
\newcommand{\ei}{\end{itemize}}
\newcommand{\ben}{\begin{enumerate}}
\newcommand{\een}{\end{enumerate}}
\newcommand{\be}{\begin{equation}}
\newcommand{\ee}{\end{equation}}
\newcommand{\bea}{\begin{eqnarray}} 
\newcommand{\eea}{\end{eqnarray}}
\newcommand{\ba}{\begin{align}} 
\newcommand{\ea}{\end{align}}
\newcommand{\bse}{\begin{subequations}} 
\newcommand{\ese}{\end{subequations}}
\newcommand{\bc}{\begin{center}}
\newcommand{\ec}{\end{center}}
\newcommand{\bfi}{\begin{figure}}
\newcommand{\efi}{\end{figure}}
\newcommand{\ca}[2]{\caption{#1 \label{#2}}}
\newcommand{\ig}[2]{\includegraphics[#1]{#2}}
\newcommand{\bp}{\begin{proof}}                 
\newcommand{\ep}{\end{proof}}                   

\newcommand{\tbox}[1]{{\mbox{\rm \tiny #1}}}
\newcommand{\ino}{\int_\Omega}
\newcommand{\mbf}[1]{{\mathbf #1}}

\DeclareMathOperator{\vol}{vol}

\newcommand{\eps}{\varepsilon}
\newcommand{\pO}{{\partial\Omega}}
\newcommand{\xx}{\mbf{x}}
\newcommand{\yy}{\mbf{y}}
\newcommand{\nn}{\mbf{n}}
\newcommand{\kk}{\mbf{k}}

\newcommand{\lo}{{L^2(\Omega)}}
\newcommand{\lpo}{{L^2(\pO)}}

\newcommand{\pois}{{\cal K}}                    
\newcommand{\si}{\sigma}                 
\newcommand{\remove}[1]{}
\newcommand{\di}[1]{{d(#1,\sigma)}}  
\newcommand{\Rn}{\mathbb{R}^n}
\newcommand{\cht}{C_\Omega} 
\newcommand{\cmp}{C_\tbox{MP}}         
\newcommand{\cb}{C_\tbox{B}}         
\newcommand{\brmk}{\begin{remark}}
\newcommand{\ermk}{\end{remark}}
\newtheorem{remark}{Remark}[section]

\newtheorem{cond}{Condition}

\newcommand\ahdelete[1]{}
\newcommand\Omegab{\partial \Omega}
\newcommand\Omegac{\overline{ \Omega}}
\newcommand\RR{\mathbb{R}}

\title{Boundary quasi-orthogonality and sharp inclusion bounds for large Dirichlet eigenvalues}

\author{A. H. Barnett and Andrew Hassell}
\date{\today}

\begin{document}
\maketitle

\begin{abstract}
We study eigenfunctions $\phi_j$ and eigenvalues $E_j$ of
the Dirichlet Laplacian on a bounded domain $\Omega\subset\RR^n$ with
piecewise smooth boundary.
We bound the distance between an arbitrary parameter $E > 0$ and the spectrum $\{ E_j \}$ 
in terms of the boundary $L^2$-norm
of a normalized trial solution $u$ of the Helmholtz
equation $(\Delta + E)u = 0$.
We also bound the $L^2$-norm of the error of this trial
solution from an eigenfunction.
Both of these results are sharp up to constants,
hold for all $E$ greater than a small constant, and
improve upon the best-known bounds
of Moler--Payne by a factor of the wavenumber $\sqrt{E}$.
One application is to the solution of
eigenvalue problems at high frequency, via, for example,
the method of particular solutions.
In the case of planar, strictly star-shaped domains we give an
inclusion bound where the constant is also sharp.
We give explicit constants in the theorems,
and show a numerical example where an eigenvalue
around the 2500th
is computed to 14 digits of relative accuracy.
%
The proof makes use of a new quasi-orthogonality property of the
boundary normal derivatives of the eigenmodes (Theorem~\ref{t:qow}
below), of interest in its own right.
Namely, the operator norm of the
sum of rank 1 operators $\partial_n \phi_j \langle \partial_n \phi_j, \cdot \rangle$
over all $E_j$ in a spectral window of width $\sqrt{E}$ --- a sum with
about $E^{(n-1)/2}$ terms --- is at most a constant factor
(independent of $E$) larger than the operator norm of any one
individual term.
\end{abstract}

\section{Introduction and main results}
\label{s:i}

\bfi 
\mbox{\raisebox{0.9in}{\parbox{3.6in}{\ig{width=3.5in}{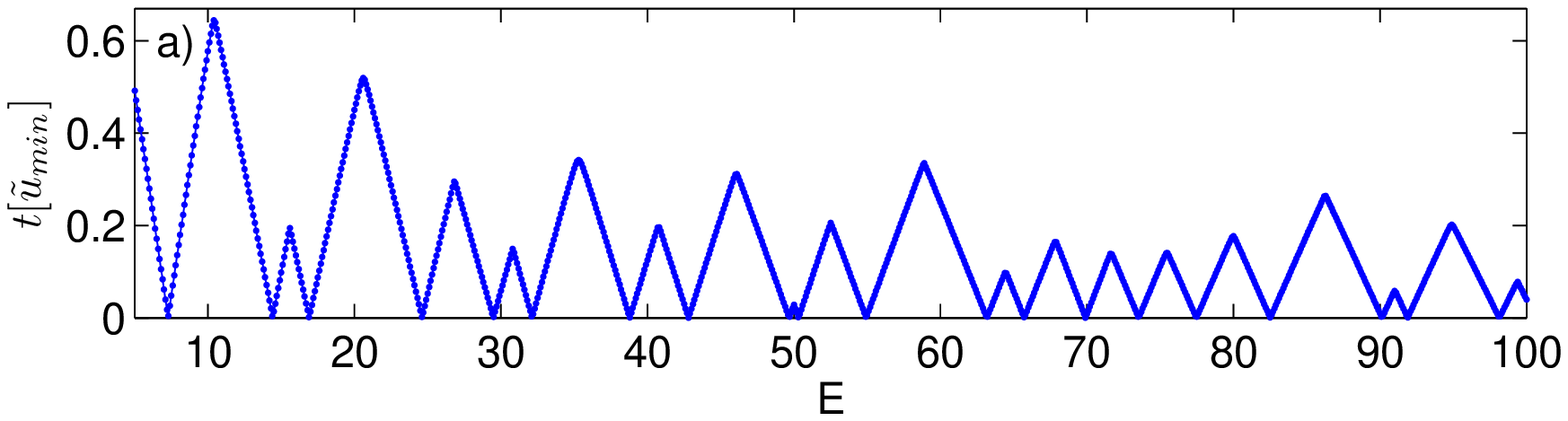}\\
\ig{width=3.5in}{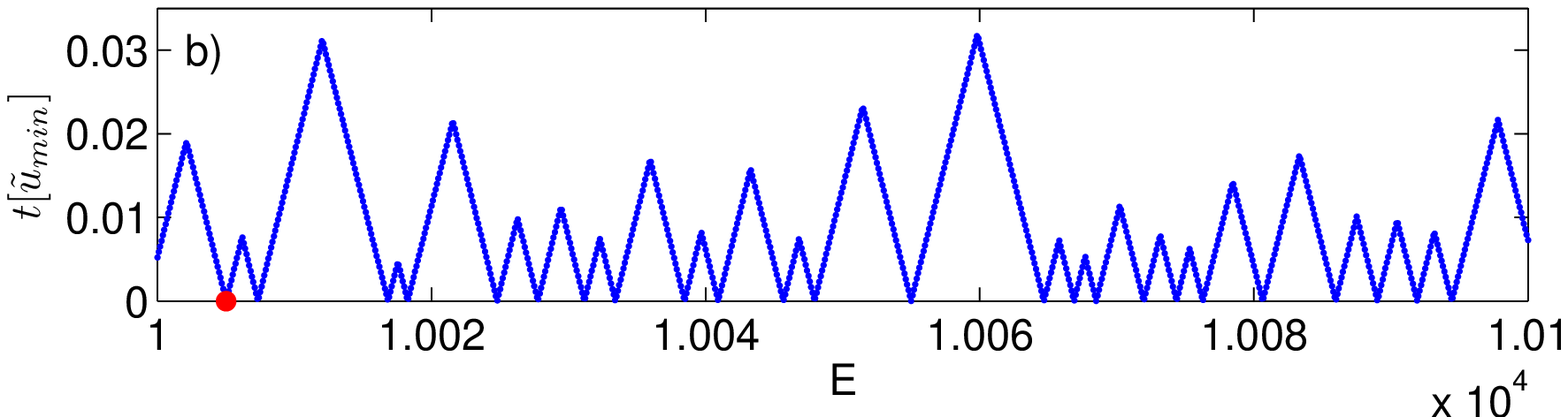}}}
\hspace{-3ex}
\ig{width=1.8in}{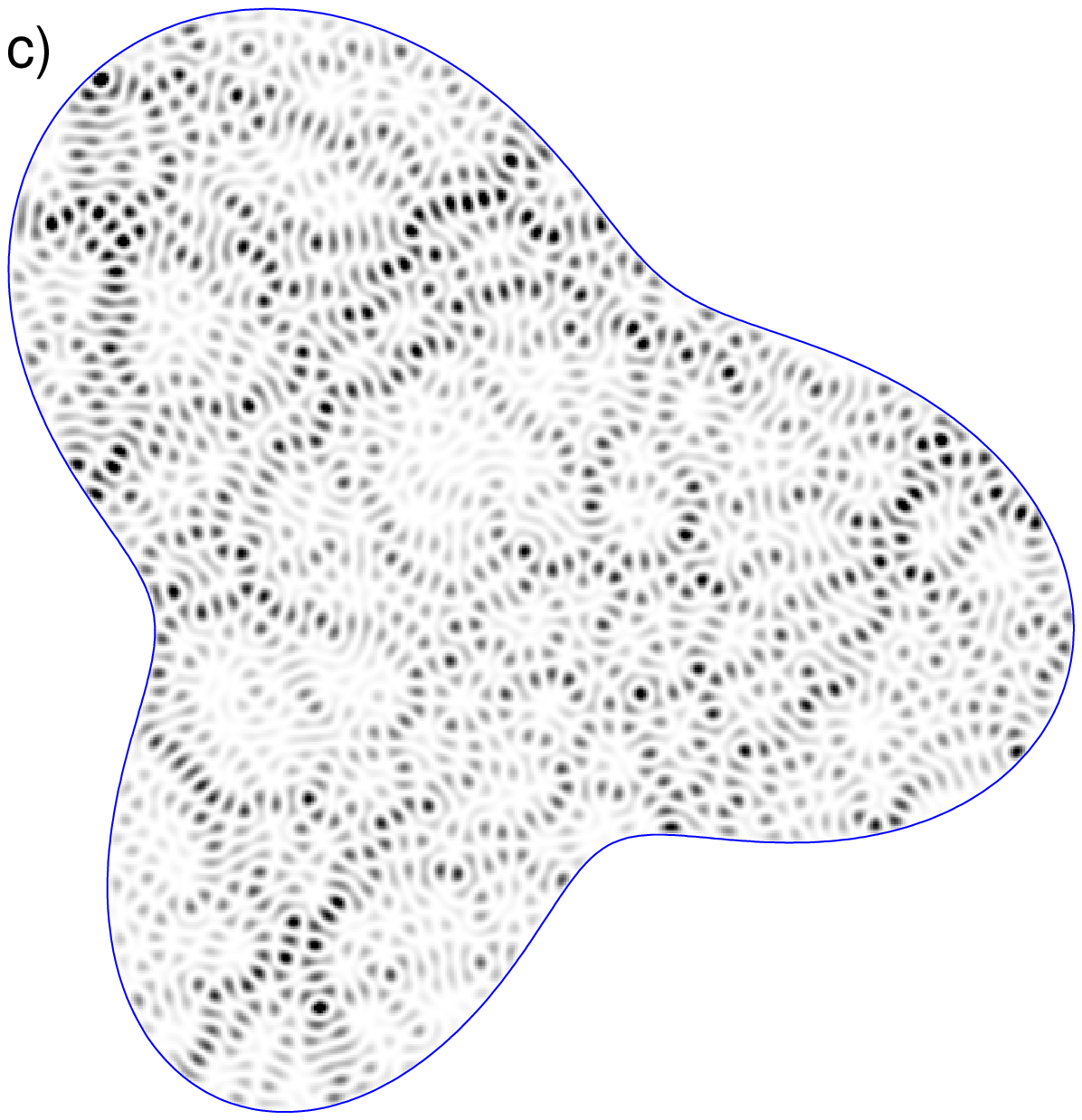}}
\ca{Tension $t[\tilde{u}_\tbox{min}]$ versus energy $E$
for the domain shown on the right.
$\tilde{u}_\tbox{min}$ is the optimal trial Helmholtz solution lying in the
span of a numerical basis set (see Section~\ref{s:num}).
a) Low frequency, showing the minima corresponding to the lowest
20 Dirichlet eigenvalues.
b) Medium-high frequency, showing a similar interval
starting at eigenvalue number $j\approx 2552$;
note the new vertical scale.
c) Density plot of eigenfunction $\phi_j \approx \tilde{u}_\tbox{min}$
corresponding to the eigenvalue $E_j=10005.02135797\cdots$
shown by the dot in b) (black indicates large values of
$|\phi_j|^2$, white zero).
}{f:t}
\efi

The computation of eigenvalues and eigenmodes of
Euclidean domains is a classical problem
(in two dimensions this is the `drum problem', reviewed in \cite{KS,tref06})
with a wealth of applications
to engineering and physics, including acoustic, electromagnetic
and optical cavity and resonator design,
micro-lasers \cite{hakan05}, and data analysis \cite{saito}.
It also has continued interest in mathematical community
in the areas of quantum chaos \cite{zencyc,que} and spectral geometry
\cite{gww}.
Let $\phi_j$ be a sequence of orthonormal eigenfunctions and $E_j$ the
respective eigenvalues ($0<E_1<E_2\le E_3\le \cdots$ counting
multiplicities) of $-\Delta$, where
$\Delta:=\sum_{m=1}^n \partial^2/\partial x_m^2$
is the Laplacian in a bounded domain $\Omega\in\Rn$, $n\ge 2$,
with Dirichlet boundary condition. That is, $\phi_j$ satisfies
\bea
(\Delta + E_j) \phi_j &=& 0 \qquad \mbox{in } \Omega
\label{e:evp}
\\
\phi_j&=&0 \qquad \mbox{on } \pO 
\\
\| \phi_j \|_{L^2(\Omega)} &=& 1.
\label{e:bc}
\eea
We will call the spectrum $\si:=\{E_j\}_{j=1}^\infty$.
Many of the applications mentioned 
demand high frequencies, that is,
mode numbers $j$ from $10^2$ to as high as $10^6$.
Efficient solution of the problem thus requires specialized numerical
approaches that scale with wavenumber
better than conventional discretization methods.

The goal of this paper is to bound the errors of approximate
eigenvalues and eigenfunctions computed using trial functions 
that satisfy exactly the homogeneous Helmholtz equation in $\Omega$.
As we will review below, such computational methods have proven very powerful.
Recently one of the authors \cite{incl} improved upon the
classical eigenvalue bound of Moler--Payne
\cite{molerpayne} by a factor of the wavenumber;
however, this result has limited utility since it applies only
to Helmholtz parameters lying in neighborhoods of $\sigma$ of unknown size.
In the present paper we go well beyond this result
by giving new theorems, which
i) hold for {\em all} Helmholtz parameters (greater than an $O(1)$ constant),
ii) retain the improved high-frequency
asymptotic behavior of \cite{incl} and show that this behavior
is sharp,
and iii) improve upon the best-known eigenfunction estimates,
again by a factor of the wavenumber.
To achieve
this we make use of a new form of quasi-orthogonality of the eigenfunctions
on the boundary, Theorem~\ref{t:qow}, of independent interest.

Before presenting our results, we need to review some known inclusion bounds
and their importance for applications.
Given an energy parameter\footnote{The Helmholtz parameter $E$ may be
interpreted as energy, or as the square of frequency, depending on the application.}
$E>0$,
let $u$ be a non-trivial solution to the homogeneous Helmholtz equation
$(\Delta+E)u=0$ in $\Omega$
with no imposed boundary condition, and define
its boundary error norm (or `tension')
\be
t[u] := \frac{\|u\|_\lpo}{\|u\|_\lo}
~.
\label{e:t}
\ee
Clearly, $t[u]=0$ implies that $E$ is an eigenvalue.
It is reasonable to expect that if $t[u]$ is small for some Helmholtz solution $u$, then $E$ is close to an eigenvalue. 
Moler--Payne \cite{molerpayne} (building upon \cite{fhm})
quantified this:
there is a constant $\cmp$ depending only on the domain, such that
\be
\di{E} \;\le\; \cmp E\, t[u]
~,
\label{e:mp}
\ee
where $\di{E}:=\min_{j}|E_j-E|$ denotes the distance of $E$ from the spectrum.

An important application is to solving \eqref{e:evp}-\eqref{e:bc}
via global approximation methods,
including the method of particular solutions (MPS) \cite{mps,incl}.
One writes a trial eigenmode
$u = \sum_{n=1}^N c_n \xi_n$ 
via basis functions $\xi_n$
which are closed-form Helmholtz solutions in $\Omega$ but which need not satisfy
any particular boundary condition.
By adjusting the coefficients $c_n$ (via a
generalized eigenvalue \cite{incl} or singular value problem \cite{gsvd})
one may minimize $t[u]$ at fixed $E$;
by repeating this
in a search for $E$ values where the minimum $t[u]$ is very small,
as illustrated by Fig.~\ref{f:t}a and b,
one may then locate approximate
eigenvalues whose error is bounded above by \eqref{e:mp}.
(This is sometimes called the method of {\em a priori-a posteriori}
inequalities \cite[Sec.~16]{KS}.)

Due to the work of Betcke--Trefethen \cite{mps} and others, such
methods have enjoyed a recent revival, at least in $n=2$,
due to their high (often spectral) accuracy
and their efficiency
at high frequency when
compared to direct discretization methods such as finite elements.
For example, in various domains, 14 digits may be achieved in double precision
arithmetic \cite{mps},
and with an MPS variant known as the
scaling method, tens of thousands of
eigenmodes as high as $j\sim 10^6$ have been computed \cite{que,mush}.
(There are also successful variants \cite{descloux,driscoll} by
Descloux--Tolley, Driscoll, and others,
in which subdomains are used, which we will not pursue here.)

If we instead interpret $u$ as the solution error for an interior
Helmholtz boundary-value problem (solved, for instance,
via MPS or boundary integral methods), then \eqref{e:mp}
states that the interior error is controlled by the boundary error;
this aids the numerical analysis of such problems \cite{litrefftz,mfs}.
Similar estimates (which, however, rely on impedance boundary conditions)
enable the analysis of least-squares non-polynomial
finite element methods \cite[Thm 3.1]{monkwang}.
Improving such estimates could thus be of general benefit
for the numerical solution of Helmholtz problems.
%

 
Recently one of the authors \cite{incl}
observed numerical evidence that \eqref{e:mp} is not sharp
for large $E$,
and showed that there is a constant $\cb$ depending only on $\Omega$,
such that, for each $\eps>0$,
\be
\di{E} \; \le \; \cb(1+\eps) \sqrt{E}\,t[u]
\label{e:incl}
\ee
holds whenever $E$ lies in some open, possibly disconnected, subset
of the real axis containing $\si$.
This is an improvement over \eqref{e:mp} by a factor of the wavenumber
$\sqrt{E}$, which
in problems of interest can be as high as $10^3$.
However, since the proof relied on analytic perturbation in the parameter $E$,
there was no knowledge about the {\em size} of this ($\eps$-dependent) subset,
hence
no way to know in a given practical situation whether the error bound holds.
The point of the present work is then
to remedy this problem
by removing any restriction to an unknown subset,
and also to extend the $\sqrt{E}$ improvement to bounds on approximate
eigenfunctions.

We assume the domain $\Omega\subset\Rn$ has unit area (or volume for $n>2$),
and obeys the following rather weak geometric condition.

\begin{cond} 
The domain $\Omega\subset\Rn$ is bounded, with piecewise smooth boundary
in the sense of Zelditch--Zworski \cite{zzw}.
This means that  $\Omega$ is given by an intersection 
$$
\Omega = \bigcap_{i=1}^N \{ \xx \mid f_i(\xx) > 0 \},
$$
where the $f_i$ 
are smooth functions defined on a neighborhood of $\overline{\Omega}$ such that 
\begin{itemize}
\item $\nabla f_i \neq 0$ on the set $\{ f_i = 0 \}$, 
\item $\{ f_i = f_j = 0 \}$ is an embedded submanifold of $\RR^n$, $1 \leq i < j \leq N$,  and
\item $\Omega$ is locally Lipschitz, i.e. for any boundary point $\xx_0 \in \Omegab$, there is a Euclidean coordinate system $z_1, \dots, z_n$ and a Lipschitz function $k$ of $n-1$ variables such that in some neighborhood of $\xx_0$,
we have 
\begin{equation}
\Omegab = \{ z_n = k(z_1, \dots, z_{n-1}) \}.
\label{e:Lip}\end{equation}
\end{itemize}
\label{c:a}
\end{cond} %

Our main result on eigenvalue inclusion is the following.

\begin{theorem} 
Let $\Omega\subset\Rn$ be a domain satisfying condition \ref{c:a}.
Then there are constants $C,c$ depending only on $\Omega$,
such that the following holds.
Let $E>1$ and suppose $u$ is a non-trivial solution of
$(\Delta+E)u=0$ in $C^\infty(\Omega)$,
with $t[u] := \|u\|_\lpo/\|u\|_\lo$.
Then,
\be
\di{E} \; \le\; C \sqrt{E}\, t[u]
~,
\label{e:b}
\ee
and for the normalized Helmholtz solution $u_\tbox{min}$ minimizing $t[u]$ at the
given $E$,
\be
c\sqrt{E}\, t[u_\tbox{min}] \; \le\; \di{E}
\; \le \; C \sqrt{E} t[u_\tbox{min}].
\label{e:ulbnds}
\ee
\label{t:b}
\end{theorem} 

\brmk
The estimate \eqref{e:ulbnds} states that \eqref{e:b} is sharp, i.e.,
using $t[u]$ alone
one cannot localize the spectrum any more tightly than this, apart from optimizing the constants $c$ and $C$.
\ermk

\brmk 
The existence of a minimizer for $t[u]$ follows from Lemma~\ref{compact}, in the case that $E$ is not a Dirichlet eigenvalue (and is trivial when $E$ is a Dirichlet eigenvalue).
The lower bound on the distance to the spectrum in \eqref{e:ulbnds}
is of use when
the numerical scheme is known to produce a good approximation to
$u_\tbox{min}$.
%
%
\ermk

We will also prove the following corresponding bound on the error of the
trial eigenfunction $u$,
which improves by a factor $\sqrt{E}$ the previous best known result
(Moler--Payne~\cite[Thm.~2]{molerpayne}).
\begin{theorem} 
Let $\Omega$ be as in Theorem~\ref{t:b}. Then there is a constant
$C$ depending only on $\Omega$, such that the following holds.
Let $E>1$, let $E_j$ be the eigenvalue nearest to $E$,
   and let $E_k$ the next nearest distinct eigenvalue.
Suppose $u$ is a solution of
$(\Delta+E)u=0$ in $C^\infty(\Omega)$ with $\|u\|_\lo=1$,
and let $\hat{u}_j$ be the projection of $u$ onto the $E_j$ eigenspace.
Then,
\be
\|u - \hat{u}_j\|_\lo \; \le \; C\frac{\sqrt{E}\,t[u]}{|E-E_k|}
~.
\label{e:e}
\ee
\label{t:e}
\end{theorem} 

\brmk
The left-hand side above is equal to
$\sin\theta$, where $\theta$ is the subspace angle
between $u$ and the $E_j$ eigenspace (this viewpoint is
elaborated in \cite[Sec.~6]{mps}).
For example, when $E_j$ is a simple eigenvalue, we may
write $\|u-\phi_j\|_\lo = 2\sin(\theta/2)$.
\ermk

\brmk
This result is also sharp, in a certain sense: see Remark~\ref{efnopt}.
\ermk

\

To conclude the introduction, we present some key ingredients
of the proofs.
Define the boundary functions of the eigenmodes by
\be
\psi_j(s) := \partial_n \phi_j(s) \qquad s\in\pO
\label{e:psi}
\ee
where $\partial_n = \nn\cdot\nabla$ is the usual normal derivative.
Our main tools will be two theorems stating that 
boundary functions $\psi_j$ lying close in eigenvalue are almost orthogonal.
The first
is the following new result which we prove in Section~\ref{s:win}. 

\begin{theorem}[spectral window quasi-orthogonality]
Let $\Omega\subset\Rn$ be a domain satisfying Condition \ref{c:a}.
There exists a constant $\cht$ depending only on $\Omega$ such that
the operator norm bound
\be
\Bigl\|\sum_{|E_j-E|\le E^{1/2}} \!\!
\psi_j\langle\psi_j,\cdot\rangle\Bigr\|_{L^2(\Omegab) \to L^2(\Omegab)}
\;\le\; \cht E
\ee
holds for all $E \geq 1$. (Here, $\langle \cdot, \cdot \rangle$ denotes the inner product in $L^2(\Omegab)$.)
\label{t:qow}
\end{theorem}

\brmk
By Weyl's Law \cite[Ch. 11]{garab} there are
$O(E^{(n-1)/2})$ terms in the above sum.
Since each term already has norm $\geq cE$ \cite{rellich, hasselltao},
the theorem expresses
essentially complete mutual orthogonality, up to a constant.
Only the scaling of the window width with $E$ is important:
the theorem also holds for a window $|E_j-E|\le c E^{1/2}$ for any fixed $c$
($\cht$ will then depend on $c$ as well as $\Omega$).
On the other hand, one could not expect it to hold over a
spectral window of width $O(E^\beta)$ for $\beta>1/2$, since
the boundary functions are approximately band-limited to spatial
wavenumber $E^{1/2}$ and thus no more than $O(E^{(n-1)/2})$
of them could be orthogonal on the boundary.
\ermk

The second result is a pairwise estimate on the inner product of boundary functions lying close in eigenvalue, with respect to a special inner product: 
(Here, $\xx(s)$ refers to the location of boundary
point $s$ relative to a fixed origin, which may or may not be inside $\Omega$.)
\begin{theorem}[pairwise quasi-orthogonality]  
Let $\Omega\subset\mathbb{R}^n$ 
be a bounded Lipschitz domain, and let $S:=
\frac{1}{2} \sup_{\xx\in\Omega} \|\xx\|$.
Then, for all $i, j\ge 1$,
\be
\Big| \int_{\pO} (\xx(s)\cdot\nn(s))\, \psi_i(s) \psi_j(s) \, ds
- 2E_i \delta_{ij} \Big|
\;\leq\;
S^2 (E_i-E_j)^2
\label{e:qo}
\ee
\label{t:qo}
\end{theorem}

\brmk
This theorem was proved by the first-named author in \cite[Appendix~B]{que}.
It  may be viewed as an off-diagonal generalization of a theorem
of Rellich \cite{rellich} which gives the $i=j$ case.
The boundary weight $\xx\cdot\nn$ (also known as the Morawetz
multiplier)
is the only one known that gives quadratic growth away
the diagonal yet also gives non-zero diagonal elements.
\ermk

Note that neither of the above quasi-orthogonality theorems implies the other.
We also note that B\"{a}cker et al.\
derived a completeness property of the boundary functions in
a (smoothed) spectral window \cite[Eq.~(53)]{backerbdry},
that is closely related to Theorem~\ref{t:qow}.

After proving Theorem~\ref{t:qow}, 
we combine it with a boundary operator defined in Section~\ref{s:A}
to prove the main theorems, in Section~\ref{s:main}.
In Section~\ref{s:star} we state and prove a variant of Theorem~\ref{t:b}
for strictly star-shaped planar domains, which has an optimal
constant $C$.
This builds on Theorem~\ref{t:qo} combined with the
Cotlar-Stein lemma (see Lemma~\ref{l:ao}).
In the main Theorems~\ref{t:b}, \ref{t:e} and \ref{t:bs},
the domain-dependent constants are not explicit;
we discuss their explicit values in Section~\ref{s:size}.
We present a high-accuracy
numerical example using the MPS, and sketch some of
the implementation aspects, in Section~\ref{s:num}.
Finally, we conclude in Section~\ref{s:conc}.


\section{Quasi-orthogonality in an eigenvalue window}
\label{s:win}

Here we prove Theorem~\ref{t:qow} using a ``$TT^\ast$ argument''.
We need the fact that the 
upper bound $\| \psi_j  \|_{L^2(\Omegab)}^2 \leq C E_j$ on eigenmode normal derivatives, proved for example in \cite{hasselltao},
generalizes to quasimodes living in an $O(E^{1/2})$ spectral window.
The proof is almost the same as in \cite{hasselltao}. 

\begin{lemma} 
Let $\Omega\subset\Rn$ satisfy Condition~\ref{c:a}.
Let $E>1$, and let
\be
\phi:=\sum_{|E_j-E|\le E^{1/2}} \!\! c_j \phi_j
\label{e:co}
\ee
with real coefficients $c_j$,  and $\sum_j c_j^2=\|\phi\|^2_\lo=1$.
Then,
\be
\|\partial_n \phi\|^2_\lpo \;\le\; \cht E
\label{e:qht}
\ee
where the constant $\cht$ depends only on $\Omega$.
\label{l:qht}
\end{lemma} 

\bp
To prove this we need the following lemma, proved in Appendix~\ref{a:outgoing},
stating that for any piecewise smooth domain (in the sense of Condition~\ref{c:a}) there is a smooth vector field that is outgoing at each boundary point. 

\begin{lemma} Let $\Omega$ satisfy Condition~\ref{c:a}. Then there exists a smooth vector field $\mbf{a}$, defined on a neighborhood of $\Omegac$, such that 
\begin{equation}
\mbf{a} \cdot \nn \geq 1
\label{outgoing}\end{equation}
 almost everywhere on $\Omegab$. 
\label{l:outgoing}\end{lemma}

The main tool for proving Lemma~\ref{l:qht} is the identity
\be
\int_\pO (D\phi)\partial_n\phi =
-\ino \phi[\Delta,D]\phi + \ino (D\phi) (\Delta+E)\phi
-\ino\phi D(\Delta+E)\phi
\label{e:com}
\ee
for any first order differential operator $D$, which follows from\footnote{The computation, involving a total of three derivatives, is justified for our class of domains, since Dirichlet eigenfunctions are in $H^{3/2}(\Omega)$ for any Lipschitz $\Omega$; see \cite{JK}, Theorem B, p164. Rellich-type computations are also justified on Lipschitz domains in \cite{Ancona}.} Green's 2nd identity, the definition of the commutator,
and $\phi|_\pO=0$.
Choosing $D:=\mbf{a}\cdot\nabla$, where $\mbf{a}$ is as in Lemma~\ref{l:outgoing},
we notice that the left-hand side of \eqref{e:com}
bounds the left-hand side of \eqref{e:qht}, since
\be
\int_\pO \psi_j^2 \; \le\; \int_\pO (\mbf{a}\cdot\nn) \psi_j^2
\ee
by Condition~\ref{c:a}.
We may now bound each of the terms on the right-hand side of \eqref{e:com}.
Defining $C_a = \sup_{\xx\in\Omega} |\mbf{a}(\xx)|$, we have
\be
\|D\phi\|_\lo^2 \le 
C_a^2 \ino \|\nabla \phi\|^2 = -C_a^2 \ino \phi \Delta \phi
= C_a^2\sum_{j} |c_j|^2 E_j \le C_a^2 F
\label{e:Dp}
\ee
where $F:=E+E^{1/2}$ is the upper end of the window. Similarly,
\bea
\|D(\Delta+E)\phi\|_\lo^2 &\le&
C_a^2\int_\Omega\|\nabla(\Delta+E)\phi\|^2
=
C_a^2
\sum_{ij}c_i c_j \int_\Omega (\Delta+E)\phi_i (-\Delta) (\Delta+E)\phi_j
\nonumber \\
&=&
C_a^2 \sum_j c_j^2 E_j (E-E_j)^2
\le C_a^2 EF
~.
\label{e:DHp}
\eea
Using Cauchy-Schwarz,
the sum of the last two terms in \eqref{e:com} is then bounded by
$2 C_a \sqrt{EF}$.
For the first term on the right of \eqref{e:com},
we use Einstein notation
$[\Delta,D] = \partial_{ii}(a_j \partial_j \cdot) - a_j \partial_{iij}$.
After several steps, using integration by parts and $\phi|_\pO=0$, we get
\be
-\int_\Omega \phi[\Delta,D]\phi =
2 \int_\Omega (\partial_i a_j)(\partial_i\phi)\partial_j\phi
+\int_\Omega (\partial_{ii}a_j)\phi\partial_j\phi
\label{e:ein}
\ee
The constants
$C'_a := \sup_{\xx\in\Omega} \|\mathbb{A}(\xx)\|_2$ where the matrix
$\mathbb{A}\in\mathbb{R}^{n\times n}$ has entries $\partial_i a_j$,
and  $C''_a := \sup_{\xx\in\Omega, j=1,\ldots,n} |\Delta a_j(\xx)|$,
exist and are finite.
Then \eqref{e:ein} is bounded by $2C'_a F + C''_a F^{1/2}$.
Adding all bounds on terms in \eqref{e:com} we get
\be
\|\partial_n \phi\|^2_\lpo \le 2(C_a +C'_a)F + C''_a \sqrt{F}
~,
\label{e:L2bnd}
\ee
which is bounded by a constant times $E$ for $E>1$.
\ep

\emph{Proof of Theorem~\ref{t:qow}.}
Consider the coefficient vector $\mbf{c}:=\{c_j\}
\in \mathbb{R}^N$ appearing in \eqref{e:co}, where $N$ is the
number of eigenvalues (counting multiplicity) in the spectral window.
Define the linear operator $T:\mathbb{R}^N\to\lpo$ by
\be
T\mbf{c} = \sum_j c_j \psi_j
\label{e:T}
\ee
Lemma~\ref{l:qht} states that $\|T\|_{l^2\to \lpo} \le (\cht E)^{1/2}$.
Thus $\|TT^\ast\|_\lpo \le \cht E$.
But $TT^\ast$ is the operator in the statement of Theorem~\ref{t:qow},
which completes its proof.

\section{Relating tension to a boundary operator}
\label{s:A}

In this section, we show, following Barnett \cite{incl}, that the
tension $t[u]$ is related to the operator norm of a natural boundary
operator.

For $E$ a non-eigenvalue of $\Omega$, let
$\pois(E):\lpo \to \lo$ be the solution operator (Poisson kernel) for the
interior Dirichlet boundary-value problem,
\bea
(\Delta + E)u &=& 0 \qquad \mbox{ in } \Omega
\label{e:helm}
\\
u &=& f \qquad \mbox{ on } \pO
\label{e:data}
~,
\eea
that is, $u = \pois f$. (For existence and uniqueness for $L^2$ data on 
a Lipschitz boundary see for example \cite[Thm.~4.25]{mclean}.)
Since the eigenbasis is complete in $\lo$,
we may write $u = \sum_{j=1}^\infty c_j \phi_j$.
We evaluate each $c_j$ by applying Green's 2nd identity,
\be
(E-E_j)(\phi_j,u)_\lo = \ino (u\Delta \phi_j - \phi_j \Delta u)
= \int_\pO (f \psi_j - \phi_j \partial_n u) ds
~,
\label{e:G2I}
\ee
thus $c_j = \langle \psi_j, f\rangle / (E-E_j)$.
The solution operator may therefore be written as a sum of rank-1 operators,
\be
\pois(E) = \sum_{j=1}^\infty \frac{\phi_j \langle \psi_j, \cdot \rangle}{E-E_j}
~.
\label{e:Ksum}
\ee
By the definition \eqref{e:t} we have, now for any $u$
satisfying $(\Delta+E)u=0$ in $\Omega$, that $t[u]^{-1} \le \|\pois(E)\|$.
Since $\|\pois^\ast\pois\| = \|\pois\|^2$, then by defining
the boundary operator in $\lpo\to\lpo$, 
\be
A(E) := \pois(E)^\ast\pois(E)
~,
\label{e:A}
\ee
we have an estimate on the tension that will be the main tool in our
analysis,
\be
t[u]^{-2}\;\le\; \|A(E)\|
~.
\label{e:tA}
\ee
Inserting \eqref{e:Ksum} into \eqref{e:A} and using orthogonality
(or see \cite[Sec.~3.1]{incl}),
we have that $A$ also may be written as the sum of rank-1 operators,
\be
A(E) = \sum_{j=1}^\infty \frac{\psi_j\langle\psi_j,\cdot\rangle}{(E-E_j)^2}
\label{e:Asum}
\ee
This sum is conditionally convergent: the sum of the operator
norm of each term diverges.
For instance, for $n=2$, Weyl's
law \cite[Ch. 11]{garab}
states that the density of eigenvalues $E_j$ is asymptotically constant,
but since $\|\psi_j\|^2=\Omega(E_j)$
the sum of norms is logarithmically divergent;
for $n>2$ the divergence is worse. 
Despite this, we have the following, which improves
upon the results of \cite{incl}.
\begin{lemma}\label{compact}
Let $\Omega\subset \Rn$, $n\ge 2$, satisfy Condition~\ref{c:a}, and let $E>0$.
Then 
\be
\lim_{N \to \infty} \ \  \sum_{j=1}^N \frac{\psi_j\langle\psi_j,\cdot\rangle}{(E-E_j)^2}
\label{e:Asum1}
\ee
converges in the
norm operator topology. Furthermore, the limit operator $A(E)$
is compact in $\lpo$.
\end{lemma}
\bp This follows immediately
from \eqref{e:Im} in the proof of Lemma~\ref{l:tail} below, which 
shows that the tail of the sum in \eqref{e:Asum} has vanishing operator norm.
$A$ is therefore also the norm limit of a sequence of finite-rank
operators.
\ep


\section{Proof of Theorems~\ref{t:b} and \ref{t:e}}
\label{s:main}

In the previous section we related tension to the norm of 
a boundary operator which itself can be written as a sum involving
mode boundary functions.
Here we place upper bounds on $\|A(E)\|$ in order to prove Theorems~\ref{t:b}
and \ref{t:e}.
Firstly we note that when $E$ is an eigenvalue, Theorem~\ref{t:b} is
trivially satisfied, since $t[u_\tbox{min}] = 0$.
When $E$ is a non-eigenvalue, formula \eqref{e:Asum}
enables us to split up contributions from different parts
of the Dirichlet spectrum,
\be
A(E) = A_\tbox{near}(E)+A_\tbox{far}(E)+A_\tbox{tail}(E)
\label{e:split}
\ee
where
\bea
A_\tbox{near}(E) &=&
\sum_{|E_j-E|\le E^{1/2}} \frac{\psi_j\langle\psi_j,\cdot\rangle}{(E-E_j)^2}
\label{e:Anear}
\\
A_\tbox{far}(E) &=&
\sum_{E/2\le E_j\le 2E, \, |E_j-E|>E^{1/2}}
\frac{\psi_j\langle\psi_j,\cdot\rangle}{(E-E_j)^2}
\label{e:Afar}
\\
A_\tbox{tail}(E) &=&
\sum_{E_j<E/2} \frac{\psi_j\langle\psi_j,\cdot\rangle}{(E-E_j)^2}
+
\sum_{E_j>2E} \frac{\psi_j\langle\psi_j,\cdot\rangle}{(E-E_j)^2}
\label{e:Atail}
\eea
It is sufficient (due to the operator triangle inequality)
to bound the norms of these three terms independently.
We first tackle the ``far'' and ``tail'' terms.

\begin{lemma}
There is a constant $C$ dependent only on $\Omega$ such that
\be
\bigl\|A_\tbox{far}(E)\bigr\| \;\le \; C \qquad \mbox{ for all } E>1
\ee
\label{l:far}
\end{lemma}
\bp
For any $E>1$,
consider the spectral interval $I_m:=[E+mE^{1/2}, E+(m+1)E^{1/2}]$.
For any such interval lying in $[E/2,2E]$ we may
apply Theorem~\ref{t:qow},  with $E$ replaced by at most $2E$,
to bound $\|\sum_{E_j\in I_m} \psi_j\langle \psi_j,\rangle\|$
by $2\cht E$.
For terms in \eqref{e:Afar} associated with this interval,
the denominators are no less than $m^2E$. Thus
\be
\Bigl\|\sum_{E_j\in I_m}
\frac{\psi_j\langle\psi_j,\cdot\rangle}{(E-E_j)^2}\Bigr\|
\;\le\;
\frac{2 \cht}{m^2}
\label{e:mwindow}\ee
Covering $[E+E^{1/2},2E]$ by summing over $m=1,2,\ldots$ gives
a constant, since $\sum m^{-2} = \pi^2/6$.
The same argument applies for intervals covering $[E/2,E-E^{1/2}]$.
\ep

\begin{lemma}
There is a constant $C$ dependent only on $\Omega$ such that
\be
\bigl\|A_\tbox{tail}(E)\bigr\| \;\le\; C E^{-1/2}\qquad\mbox{ for all } E>1
\ee
\label{l:tail}
\end{lemma}
\bp
Consider a spectral interval $I_m:=[2^m E, 2^{m+1} E]$. We may cover
this with at most $2^{m/2 -1}E^{1/2} + 1$ windows of half-width
at most $2^{m/2}E^{1/2}$; for each of these windows Theorem~\ref{t:qow}
applies to bound $\|\sum_{E_j\in I_m} \psi_j\langle \psi_j,\rangle\|$
by $\cht 2^{m+1}E$. 
For each $E_j\in I_m$, the denominator is no smaller than $(2^{m-1}E)^2$.
Thus
\be
\Bigl\|\sum_{E_j\in I_m}
\frac{\psi_j\langle\psi_j,\cdot\rangle}{(E-E_j)^2}\Bigr\|
\;\le\;
(2^{m/2 -1}E^{1/2} + 1) \frac{\cht 2^{m+1} E}{(2^{m-1}E)^2}
= \cht (2^{-m/2 -2} E^{-1/2} + 2^{-m+1} E^{-1})
\label{e:Im}
\ee
The infinite sum over $m=1,2,\ldots$ gives
\be
\Bigl\|\sum_{E_j>2E} \frac{\psi_j\langle\psi_j,\cdot\rangle}{(E-E_j)^2}\Bigr\|
\le
\cht\left(\frac{E^{-1/2}}{4(\sqrt{2}-1)} + 2E^{-1}\right) \le C E^{-1/2}
\; \mbox{ for all } E>1.
\label{e:toptail}
\ee
We treat the interval $(0,E/2)$ similarly, using a sequence
of intervals $J_m:=[2^{-m-1}E,2^{-m}E]$. Each such
interval may be covered by at most $2^{-(m+3)/2} E^{1/2} +1$ windows
of half-width $2^{-(m+1)/2}E^{1/2}$. For each $E_j\in J_m$,
the denominator is no smaller than $E^2/4$. In a similar manner as before, the
operator norm of the partial sum associated with $J_m$ is then
$O(2^{-m}E^{-1/2})$, thus the infinite sum over $m$ is $O(E^{-1/2})$.
Note that Theorem~\ref{t:qow} does not apply for $E<1$, but that there
are $O(1)$ such $E_j$ values and each contributes $O(E^{-1})$.
This proves the Lemma.
\ep

\emph{Proof of Theorem~\ref{t:b}.} 
Examining the ``near'' term \eqref{e:Anear}, we use
Theorem~\ref{t:qow} on the sum of numerators, and get a bound by taking the
minimum denominator,
\be
\bigl\|A_\tbox{near}(E)\bigr\|
\;\le\;
\frac{\cht E}{\di{E}^2} \qquad \mbox{ for all } E>1
\label{e:Anearbnd}
\ee
Using this and the above Lemmas to sum the terms in \eqref{e:split} 
gives
\be
\bigl\|A(E)\bigr\|
\;\le\;
\frac{\cht E}{\di{E}^2} + C \qquad \mbox{ for all } E>1
\label{e:Abnd2}
\ee
From Lemma~\ref{l:dist}, an upper bound on the distance to the spectrum,
we see that the second term is bounded by
at most a constant times the first, so may be absorbed into it to give
\be
\bigl\|A(E)\bigr\|
\;\le\;
\frac{C E}{\di{E}^2} \qquad \mbox{ for all } E>1
\label{e:Abnd}
\ee
Combining this with \eqref{e:tA} proves \eqref{e:b}, hence also the
second inequality in \eqref{e:ulbnds}.
The first inequality in  \eqref{e:ulbnds}
simply follows from the fact that, since $A$ is a sum of positive operators,
\be
t[u_\tbox{min}]^{-2} = \bigl\|A(E)\bigr\|
\ge
\Bigl\| \frac{\psi_j\langle\psi_j,\cdot\rangle}{(E-E_j)^2}\Bigr\|
=
\frac{\|\psi_j\|^2}{\di{E}^2}
~,
\label{e:tumin}
\ee
where $E_j$ is the eigenvalue closest to $E$.
Using the lower bound $\|\psi_j\|^2\ge c E_j$ from \cite{hasselltao}
this becomes
\be
\di{E} \ge c \sqrt{E_j} t[u_\tbox{min}]
~.
\label{e:dtumin}
\ee
With a change of constant,
$E_j$ may be replaced here by
$E$ to give the first inequality in  \eqref{e:ulbnds},
since Lemma~\ref{l:dist} insures that $E_j$ is relatively close to $E$.
(The lemma is not useful for $E$ less than some constant and $E_j<E$, but then
the ratio $E/E_j$ is still bounded by a constant because $E_j\ge E_1$).

\emph{Proof of Theorem~\ref{t:e}.} 
We next prove the eigenfunction error bound \eqref{e:e},
first considering $E$ a non-eigenvalue.
We denote the boundary data by $U:=u|_\pO$.
From orthogonality,
then using the formula for the $c_i$ coefficients below \eqref{e:G2I},
we get,
\be
\|u-\hat{u}_j\|^2_\lo = \sum_{E_i\ne E_j} |(\phi_i,u)_\lo|^2
=\sum_{E_i\ne E_j} \frac{|\langle \psi_i, U \rangle|^2}{(E-E_i)^2}
\le \Bigl\|\sum_{E_i\neq E_j}
\frac{\psi_i\langle\psi_i,\cdot\rangle}{(E-E_i)^2}\Bigr\|
\|U\|_2^2 .
\label{e:proj}
\ee
The operator in the last expression is identical to \eqref{e:Asum} except with
the $E_j$-eigenspace terms omitted.
Therefore, its norm may be bounded in the same way as that of $A(E)$,
the only difference being that 
the $\di{E}$ introduced in \eqref{e:Anearbnd} is replaced by
$\min_{E_i\neq E_j}|E-E_i| = |E-E_k|$.
Thus the bound analogous to \eqref{e:Abnd} is
\[
\bigl\|\sum_{E_i\neq E_j}
\frac{\psi_j\langle\psi_j,\cdot\rangle}{(E-E_j)^2}\bigr\|
\;\le\; \frac{CE}{(E-E_k)^2}
\qquad \mbox{ for all } E>1
~,
\]
and inserting this and $\|U\|_\lpo = t[u] \|u\|_\lo = t[u]$ into
\eqref{e:proj} gives \eqref{e:e}.
Finally, if $E$ is an eigenvalue, i.e. $E=E_j$,
the solution operator
\eqref{e:Ksum} is undefined,
since a solution to \eqref{e:helm}-\eqref{e:data} exists if and only if
$f$ is orthogonal to the normal derivative functions in the $E$-eigenspace.
This can be seen by applying Green's 2nd identity to $\phi$, any function
in the $E$-eigenspace, and $u$, giving $\langle \partial_n \phi, U \rangle=0$.
However, the solution
coefficients $c_i$ for which $E_i\neq E$ are uniquely defined
by the same formula as before.
Thus \eqref{e:proj} and the rest of the proof carries through.

\brmk\label{efnopt}
Theorem~\ref{t:e} is sharp, as can be seen in the following way:
if $u$ is such that $t[u]$ is close to $t[u_{\tbox{min}}]$ (say, less than $2 t[u_{\tbox{min}}]$), then we have, by combining \eqref{e:ulbnds} and \eqref{e:e}, 
\be
\|u - \hat{u}_j\|_\lo \; \le \; C\frac{|E-E_j|}{|E-E_k|}
~.
\label{e:e2}
\ee
Apart from the value of the constant, one cannot expect to do better than this. For example, if $E$ is midway between $E_j$ and $E_k$, then the error $\|u - \hat{u}_j\|_\lo$ cannot be expected to be better than $1/\sqrt{2}$. 
\ermk

\section{Star-shaped planar domains}
\label{s:star}

The purpose of this section is to say something stronger than Theorem~\ref{t:b}
in the special case of star-shaped domains
in $n=2$. We take weighted boundary functions
\be
\psi^{(s)}_j(s) := (\xx(s)\cdot\nn(s))\partial_n \phi_j(s)
~, \qquad s\in\pO
\ee
and our boundary inner product as
\be
\langle f,g\rangle_s := \int_\pO (\xx(s)\cdot\nn(s))^{-1} f(s) g(s) ds
\ee
hence norm $\|f\|_s := \sqrt{\langle f,f\rangle}$, and
$t_s[u]:= \|U\|_s/\|u\|_\lo$.
The significance of the weight $(\xx\cdot\nn)$ is twofold:
it is strictly positive for strictly star-shaped domains, and
also turns the inner product in \eqref{e:qo} into $\langle \psi^{(s)}_i,
\psi^{(s)}_j\rangle_s$, enabling us to benefit from pairwise
quasi-orthogonality.
The Rellich theorem $\|\psi^{(s)}_j\|_s^2 = 2E_j$ states that,
with this special weight,
there is no fluctuation in the $L^2$-norms of the boundary functions.
As shown in \cite{incl}, the function $t_s[u_\tbox{min}]$ vs
$E$ has slope $1/\|\psi^{(s)}_j\|_s^2$ in the neighborhood of $E_j$
(this arises from dominance of a single term in \eqref{e:Assum} below).
Hence these slopes are predictable
{\em independently} of the particular form of each mode $\phi_j$.
This enables us to get the following eigenvalue inclusion result
analogous to Theorem~\ref{t:b}.
\begin{theorem} 
Let $\Omega\subset\mathbb{R}^2$
be a strictly star-shaped bounded domain with piecewise
smooth boundary. Then there are constants $c_1$, $c_2$, $c_3$
depending only on $\Omega$,
such that the following holds. Let $E>1$,
and suppose $u$ is a non-trivial solution of $(\Delta+E)u=0$ in
$C^\infty(\Omega)$, with $c_2t_s[u]^2<1$. Let $F:=E+\sqrt{E}$. Then,
\be
\di{E} \; \le \; \sqrt{2F}\, t_s[u]
\frac{1 + c_1 \sqrt{F} t_s[u]}{1-c_2 t_s[u]^2}
~.
\label{e:bs}
\ee
For the Helmholtz solution $u_\tbox{min}$ minimizing $t_s[u]$ at the
given $E$,
\be
\sqrt{2(E-c_3E^{1/2})}\,t_s[u_\tbox{min}] \; \le \;
\di{E}
\label{e:lbndss}
\ee
\label{t:bs}
\end{theorem} 

\brmk
In the limit of high frequency $E\gg 1$ and small tension $t_s[u] \ll E^{-1/2}$,
the right-hand side of \eqref{e:bs} and the left hand side of \eqref{e:lbndss} are both $\sqrt{2E}(1 + o(1))t_s$.
This proves that both the power of $E$ and the constant $\sqrt{2}$ are sharp.
\ermk

\brmk
Notice that this theorem is not applicable for all $E$ since there may
be large spectral gaps where $c_2t_s[u]^2<1$ cannot be satisfied.
Due to the numerator, it becomes far from optimal when $t_s[u]$
is $O(E^{-1/2})$ or larger.
In these respects it is less general than Theorem~\ref{t:b},
even though it gives better bounds in the small tension limit.
\ermk


The main tool used in the proof of Theorem~\ref{t:bs} is the pairwise quasi-orthogonality result, Theorem~\ref{t:qo}, together with the Cotlar-Stein lemma, which we state here for the special case of self-adjoint operators:

\begin{lemma}[Cotlar-Stein~\cite{cotlar,stein,Comech}]
Let $\{T_j\}_{j\in J}$ be a countable set of bounded self-adjoint operators,
$J \subset \mathbb{N}$. Then
\[
\Bigl\|\sum_{j\in J} T_j \Bigr\| \le 
\max_{j\in J} \sum_{i\in J} \sqrt{\| T_i T_j \|}
~.
\]
\label{l:ao}
\end{lemma}

\emph{Proof of Theorem~\ref{t:bs}.}
The weighted equivalent of \eqref{e:Asum} is the operator
\be
A^{(s)}(E) = \sum_{j=1}^\infty \frac{\psi^{(s)}_j\langle\psi^{(s)}_j,\cdot\rangle_s}{(E-E_j)^2}
\label{e:Assum}
\ee
which, by analogy with \eqref{e:tA}, satisfies
\be
t_s[u]^{-2}\;\le\; \|A^{(s)}(E)\|_s
~.
\label{e:tsAs}
\ee
The lower bound \eqref{e:lbndss} follows by analogy with
\eqref{e:tumin}-\eqref{e:dtumin},
using $\|\psi^{(s)}_j\|_s^2 = 2E_j$, and $E_j\ge E-c_3 E^{1/2}$
from Lemma~\ref{l:dist}.

Using the same splitting into ``near'', ``far'', and ``tail'' parts as
in Section~\ref{s:main}, we can bound the norm of the ``near'' part
in a new way, as follows.
\begin{lemma} 
There is a constant $c_1>0$ depending only on $\Omega$ such that
\[
\|A^{(s)}_\tbox{near}(E)\|_s \le \frac{2F}{\di{E}^2} +
\frac{\sqrt{2}c_1 F}{\di{E}}
\qquad \mbox{ for all } E > 1
\]
\label{l:Anear}
\end{lemma}
The first term in this bound will arise simply from the
single term in the sum \eqref{e:Assum} with $E_j$ nearest to $E$.
The second term requires more work, as we now show.
\bp  
Let $J = \{j: |E_j-E|\le E^{1/2}\}$.
Using $T_j = \frac{\psi_j\langle\psi_j,\cdot\rangle}{(E-E_j)^2}$
in Lemma~\ref{l:ao} gives
\be
\Bigl\|\sum_{j\in J} T_j\Bigr\|_s
\;\le\;
\max_{j\in J} \frac{\|\psi^{(s)}_j\|_s^{1/2}}{|E-E_j|}
\sum_{i\in J}
\frac{ \bigl(\langle\psi^{(s)}_i,\psi^{(s)}_j \rangle_s \,
\|\psi^{(s)}_i\|_s\bigr)^{1/2}}
{|E-E_i|}
\label{e:aoapply}
\ee
Applying quasi-orthogonality (Theorem~\ref{t:qo}) for the inner product,
and $\|\psi^{(s)}_j\|_s^2=2E_j$, and separating diagonal ($i=j$) from off-diagonal terms, we get,
\be
\Bigl\|\sum_{j\in J}\frac{\psi_j\langle\psi_j,\cdot\rangle}{(E-E_j)^2} \Bigr\|
\;\le\;
\max_{j\in J} \frac{2E_j}{(E-E_j)^2} +
\sqrt{2} S \max_{j\in J} \frac{E_j^{1/4}}{|E-E_j|}
\sum_{i\in J}
\frac{E_i^{1/4}|E_i-E_j|}{|E-E_i|}
\label{e:qoao}
\ee
Here $S$ is as in Theorem~\ref{t:qo}. 
The first term is bounded by $2F/\di{E}^2$. Using
$|E_i-E_j| \le |E_i-E| + |E-E_j|$
bounds the second term by
\bea
&&\sqrt{2} S \max_{j\in J} \frac{E_j^{1/4}}{|E-E_j|}
\sum_{i\in J}E_i^{1/4}\left(1+\frac{|E-E_j|}{|E-E_i|}\right)
\nonumber \\
&\le&
\sqrt{2F} S\,|J|\, \max_{j\in J}
\Bigl( \frac{1}{|E-E_j|} + \frac{1}{\di{E}}\Bigr)
\;\le\;
\frac{2\sqrt{2F} S|J|}{\di{E}}
\label{e:term2}
\eea
Recall Weyl's law for the asymptotic density of eigenvalues,
which states that, for $n=2$ and $\vol\Omega=1$, 
\be
N(E) := \#\{j:E_j<E\} = \frac{1}{4\pi}E + R(E)
~,
\label{e:weyl}
\ee
where the remainder is $R(E)=O(\sqrt{E})$ (\cite{SafVas};
for the case of piecewise-smooth boundary see \cite[Eq.~(0.3)]{Seeley}).
Since the remainder is bounded for small $E$,
there is a constant $C_\tbox{W}$ such that $|R(E)| \le C_\tbox{W} \sqrt{E}$
for all $E>1$.
Thus $|J|$, the number of terms in the ``near'' window, is bounded by
\[
|J| \le \Bigl(\frac{1}{4\pi} + 2C_\tbox{W}\Bigr) \sqrt{F}
~.
\]
Inserting this into \eqref{e:term2} proves the Lemma,
and we may take $c_1 = 2S(1/4\pi + 2C_\tbox{W})$.
\ep

\emph{Completion of the proof of Theorem~\ref{t:bs}.}
The proofs of analogously weighted versions of
Lemmas~\ref{l:far} and \ref{l:tail} are unchanged.
So we may combine them with Lemma~\ref{l:Anear} and \eqref{e:tsAs}
to get, for some constant $c_2$,
\[
t_s[u]^{-2} \;\le\; \|A^{(s)}(E)\|_s \;\le\;
\frac{2F}{\di{E}^2} + \frac{\sqrt{2}c_1 F}{\di{E}} + c_2
\qquad \mbox{ for all } E > 1
\]
Multiplying through by $\di{E}^2$
we solve the quadratic inequality for $\di{E}$,
\[
\di{E} \; \le\;
\frac{c_1 F/\sqrt{2} + \sqrt{c_1^2 F^2/2 + 2(t_s[u]^{-2}-c_2)F}}
{t_s[u]^{-2}-c_2}
\]
Using the subadditivity of the square-root completes the
proof of \eqref{e:bs}.

\section{Discussion of explicit constants} 
\label{s:size}

For the practical application of Theorems~\ref{t:b} and \ref{t:e},
it is important to have an explicit value for the constant $C$
(from the discussion after \eqref{e:proj} we notice that
$C$ in the two theorems is the same.)

We now compute an explicit value of this $C$ that holds for all $E>1$.
Examining \eqref{e:L2bnd} we see that a choice of constant in
Lemma~\ref{l:qht}, and hence Theorem~\ref{t:qow}, that holds for
all $E>1$ is $\cht = 4(C_a+C'_a) + \sqrt{2}C''_a$.
To compute this we need sup norms of the value, and first and second
derivative, of a vector field $\mbf{a}$ as in Lemma~\ref{l:outgoing}. 
                                            
The proof of Lemma~\ref{l:outgoing} shows such a construction; the values will
depend on the size of the vectors $\mbf{a}_{x_i}$ and the choice of partition of unity used to cover $\overline{\Omega}$. The vectors $\mbf{a}_{x_i}$ will be large (order $1/\epsilon$) if $\Omega$ has corners with angles less than $\epsilon$ or greater than $2\pi - \epsilon$. 
We note that a numerical procedure for this construction
could be useful.

In some special cases,  a simpler prescription for the vector field can be given:
\begin{itemize} 
\item 
For strictly star-shaped domains in $\Rn$, we may choose 
$\mbf{a}=\xx / \inf_\pO(\xx\cdot\nn)$, which 
gives $C_a =  \sup_\pO (\xx\cdot\nn)/ \inf_\pO(\xx\cdot\nn)$,
$C'_a=1 /  \inf_\pO(\xx\cdot\nn)$, and $C''_a=0$.
\item 
For a domain with $C^2$ boundary,
let $\delta > 0$ be the largest number such that for each $\xx_0 \in \Omegab$, a ball of radius $\delta$ can be placed within $\Omega$ so as to be tangent to $\Omegab$ at $\xx_0$.
We may then
choose $\mbf{a} = (1-r/\delta)^2\nn_r$, for $r<\delta$, $\mbf{a}=\mbf{0}$
otherwise, where the coordinate $r$ is the distance from $\pO$,
and $\nn_r$ is the unit vector in the local decreasing $r$ direction.
This gives constants $C_a = 1$ and $C'_a = 2/\delta$.
$C''_a$ depends on $\delta$ and an upper bound on the rate of change
of surface curvature.
(Also note that a slight modification of the proof of Theorem~\ref{t:qow}
would allow estimation purely in terms of $C_a$ and $C'_a$, but with a doubling
of the numerical constants).
\end{itemize}

Summing the terms \eqref{e:mwindow} above and below $E$ we have that the
constant in Lemma~\ref{l:far} is $2\pi^2\cht/3$.
Similarly, using \eqref{e:toptail} and its equivalent for $(0,E/2)$
gives the constant in Lemma~\ref{l:tail} as
$\cht(\frac{1}{4(\sqrt{2}-1)}+\frac{1}{4-\sqrt{2}}+6)<7\cht$.
Summing these two constants gives
a constant $C$ in \eqref{e:Abnd2} as $14\cht$.
A choice of constant in \eqref{e:Abnd} is then $\cht + 14 \cht \max[E_1^2,C_d^2]$,
where from Appendix~\ref{a:dist} we have $C_d = 2\sqrt{E_1}$, and the max
accounts for the case $1<E\le E_1$.
Finally, the constant in \eqref{e:b} is the square-root of this,
$C = \sqrt{\cht(1 + 14\max[E_1,4]E_1)}$.

Requiring that the above estimates hold for all $E>1$
caused non-optimality in the choice of constant.
It is more sensible in high frequency applications
to use a better constant which is approached for $E\gg1$, and small
tension $t\ll1$. We now give this explicitly.
In this limit, in \eqref{e:L2bnd}, $F$ tends to $E$,
and we drop lower-order terms to get $\cht = 2(C_a+C'_a)$,
which in the star-shaped case is
\be
\cht = 2\frac{1+\sup_\pO(\xx\cdot\nn)}{\inf_\pO(\xx\cdot\nn)}
\qquad \mbox{ for } E\gg 1, \mbox{ $\Omega$ star-shaped .}
\label{e:chtsshf}
\ee
If tension is small (i.e.\ $E$ is not in a large spectral gap), the
second term in \eqref{e:Abnd} is negligible, so 
we may approximate the constant in \eqref{e:b} as
\be
C = \sqrt{\cht} \qquad \mbox{ for } E\gg 1, \;t\ll 1
~.
\label{e:chf}
\ee
\brmk
The limiting constant \eqref{e:chf} does not reflect the
limiting slopes of the graph $t[u_\tbox{min}]$ vs $E$ near eigenvalues.
These slopes are known \cite{incl} to be
$1/\|\psi_j\|^2$, which is bounded by $(2C_a' E_j)^{-1}$ \cite{hasselltao}.
\ermk


We end by discussing the constants $c_1$ and $c_2$ in Theorem~\ref{t:bs}. Constant 
$c_2$ may be estimated easily, as above, using the weighted
versions of Lemmas~\ref{l:far} and \ref{l:tail}.
In the proof of Lemma~\ref{l:Anear}, $c_1$
involves the Weyl constant $C_\tbox{W}$;
we know of no explicit estimates for $C_\tbox{W}$ in the literature
(the closest we know are estimates
of the form $|R(E)| < C \sqrt{E}\ln E$ with explicit constants
\cite{netrusov,chenhua}).
However, these constants are effectively
irrelevant for practical purposes, when $E\gg1$ and $t\ll E^{-1/2}$,
since in these limits, one may replace \eqref{e:bs}
by $\di{E}\le \sqrt{2E} t_s[u]$ and still have an error bound very close
to that given by the full expression.

\section{Numerical example} 
\label{s:num}

In Fig.~\ref{f:t}c we show a planar nonconvex
domain given by the radial function
$r(\theta) = 1 + 0.3 \cos[3(\theta+ 0.2\sin\theta)]$.
The domain is star-shaped and smooth (we will not address numerical
issues raised by corners here; see
\cite{fhm,descloux,mps,driscoll,que,timodd}.)
For high-frequency eigenvalue problems,
a convenient computational basis of Helmholtz solutions are
`method of fundamental solutions' basis functions
$\xi_n(\xx) = Y_0(\sqrt{E}|\xx-\yy_n|)$,
where $Y_0$ is the irregular Bessel function of order zero,
and
$\{\yy_n\}_{n=1}^N$ are a set of `charge points'
in $\mathbb{R}^2\setminus\overline{\Omega}$. The latter were chosen
by a displacement of the boundary parametrization
$\xx(\theta)$, $0<\theta\le2\pi$, in the imaginary direction (see \cite{mfs});
specifically $\yy_n = \xx(2\pi n/N - 0.025 i)$.

We compute the data plotted in Fig.~\ref{f:t}a, b as follows.
At each $E$,
$t[\tilde{u}_\tbox{min}]$ is given by the square-root of the smallest
generalized eigenvalue of
a generalized eigenvalue problem (GEVP) involving
$N\times N$ symmetric real dense matrices $F$ and $G$
(the basis representations of the boundary and interior norms
respectively.)
Both matrices are evaluated using $M$-point
periodic trapezoidal quadrature in $\theta$, that is,
quadrature points $\xx_m = \xx(2\pi m/M)$,
$m=1,\ldots,M$, and weights $w_m = 2\pi|\xx'(2\pi m/M)|/M$.
For instance, $F = P^\ast P$, where $P\in\mathbb{R}^{M\times N}$ has elements
\be
P_{mn} = \sqrt{w_m} \xi_n(\xx_m)
~,
\label{e:P}
\ee
and $G$ is similarly found \cite[Sec.~4.1]{incl}
using $P$ and the matrices
$P^{(1)}$ and $P^{(2)}$ whose entries are the $x_1$- and $x_2$-derivatives
of those in $P$.
Since the GEVP is numerically singular,
regularization was first performed, similarly to \cite[Sec.~6]{gsvd}, by
restricting to an orthonormal basis for the numerical column space of
$[P;P^{(1)};P^{(2)}]$ given by the
left singular vectors with singular values at least $10^{-14}$ times the
largest singular value.

For low frequencies (Fig.~\ref{f:t}a), 8-digit accuracy requires
$N=100$ basis functions and $M=200$ quadrature points. For higher frequencies
corresponding to 40 wavelengths across the domain
(Fig.~\ref{f:t}b, c), it requires
$N=400$ and $M=500$,
and the above GEVP procedure takes 3 seconds per $E$ value.%
  \footnote{All computation times are reported for a laptop with
    2GHz Intel Core Duo processor and 2GB RAM, running MATLAB 2008a
    on a linux kernel.}
Very small ($<10^{-8}$) tensions cannot be found this way,
and instead are best
approximated via the GSVD \cite{gsvd}:
the optimal tension at a given $E$ is the lowest generalized singular value
of the matrix pair $(P,Q)$, where matrix $Q$ has entries
\be
Q_{mn} = \sqrt{\frac{\xx_m\cdot \nn_m}{2E}} \sqrt{w_m}
\,\frac{\partial \xi_n}{\partial n}(\xx_m)
~,
\label{e:Q}
\ee
where $\nn_m$ is the normal at $\xx_m$, and regularization as before.
Note that $G = Q^\ast Q$ well
approximates the interior norm in the subspace with zero Dirichlet data,
due to the Rellich formula (case $i=j$ of Theorem~\ref{t:qo}).

Any single-variable function minimization algorithm may then be used to
search for a local minimum of $t[\tilde{u}_\tbox{min}]$ vs $E$;
we prefer iterated fitting of
a parabola to $t[\tilde{u}_\tbox{min}]^2$ at three nearby $E$ values,
which converges typically in 5 iterations.
Using this with the GSVD
(with $N=500$, $M=700$, i.e.\ 6 points per wavelength on
$\pO$, and taking 8 seconds per iteration),
we find the tension
\be
t[\tilde{u}_\tbox{min}] = 2.2\times 10^{-12}
\qquad \mbox{at } E = 10005.0213579739
~.
\label{e:dp}
\ee
This is shown by the dot in Fig.~\ref{f:t}b.
The GSVD right singular vector gives the basis coefficients
of the corresponding trial function $\tilde{u}_\tbox{min}$, which is
plotted in Fig.~\ref{f:t}c (this took 34 seconds to evaluate on a square grid
of size 0.005, i.e.\ $1.3\times 10^5$ points.)

Armed with datapoint \eqref{e:dp}, what
can we deduce about Dirichlet eigenpairs of $\Omega$
using our new theorems, and how much better are they than
previous results?
The constant in the Moler--Payne bound
\eqref{e:mp} is $\cmp = q_1^{-1/2}$ where $q_1$ is the
lowest eigenvalue of a Stekloff eigenproblem on $\Omega$
\cite[(2.11)]{KSstek}.
Since $\Omega$ is star-shaped, the bound
$q_1 \ge E_1^{1/2}\inf_\pO(\xx\cdot\nn)/2\sup_\pO(\xx\cdot\nn)$
from \cite[Table I]{KSstek} applies,
giving $\cmp = 1.31$ as a valid choice.
Thus \eqref{e:mp} states that there is an
eigenvalue $E_j$ a distance no more than $2.9\times10^{-8}$ from the above $E$.
On the other hand, \eqref{e:chf} and \eqref{e:chtsshf} give
the constant in Theorem~\ref{t:b} as $C
= 2.9$.
Applying the theorem gives a distance from the spectrum of no more than
$6.3\times 10^{-10}$.
Taking the small-tension limit of the star-shaped planar result \eqref{e:bs},
and recomputing the weighted tension $t_s[\tilde{u}_\tbox{min}]$
from Section~\ref{s:star},
we get an even smaller distance of
$3.5\times 10^{-10}$, that is, an error of $\pm3$ in
the last digit of \eqref{e:dp}.
The latter is a 80-fold improvement over Moler--Payne.
(Also see \cite{incl}
for an example at higher frequency with 3 digits of improvement).

How good an approximation is $\tilde{u}_\tbox{min}$ to the
eigenfunction $\phi_j$?
Using the observation that the next nearest eigenvalue is $E_k
=10007.339\cdots$,
the eigenfunction bound of
Moler--Payne \cite{molerpayne} gives an $L^2$-error of $1.2\times 10^{-8}$.
With the same data, using the constant $C$ above,
Theorem~\ref{t:e} gives an $L^2$-error of $2.7\times 10^{-10}$,
a 50-fold improvement over that achievable with previously known theorems.

\section{Conclusions}
\label{s:conc}

We have improved, by a factor of the wavenumber,
the Moler--Payne bounds on Dirichlet eigenvalues and eigenfunctions
which have been the standard for the last 40 years.
This makes rigorous the conjectures
based on numerical observations in \cite{incl}.
We expect this to be useful
since high-frequency wave and eigenvalue calculations are
finding more applications in recent years.
Of independent interest is a new quasi-orthogonality result in a spectral window
(Theorem~\ref{t:qow}).

For numerical utility,
throughout we have been explicit with constants, and
have specified a lower bound on $E$ for which the estimates hold
(this being stronger than merely a `big-O' asymptotic estimate).
For star-shaped domains we strengthened the inclusion bounds
(Theorem~\ref{t:bs}), achieving a sharp power of $E$ and sharp constant,
in the limit of small tension, when tension is weighted by
a special geometric function.
This weight allowed pairwise quasi-orthogonality to be used, but
since an upper bound for the number of eigenvalues
in a $\sqrt{E}$ window is needed,
this is only useful for $n=2$ (planar domains).

We applied our theorems to a numerical example,
enabling close to 14 digits accuracy in a high-lying eigenvalue,
and 10 digits in the eigenfunction. Both are two digits beyond
what could be claimed with previously-known theorems.

Our estimate $C \sqrt{E} t[u]$ on the distance to the spectrum is
sharp (up to constants) if the tension $t[u]$ (or $t_s[u]$) is comparable to $t[u_{\tbox{min}}]$ ($t_s[u_{\tbox{min}}]$). 
However, numerically, one generally has access to other properties
of $u$ (e.g. its normal derivative)
which can give more detailed information  about the spectrum.
For example, the powerful `scaling method' \cite{v+s,que}
 is able to locate many eigenvalues using an operator computed
at a single energy $E$. In another direction, 
 Still \cite[Thm.~4]{still} obtains improved inclusion bounds when the
 approximate eigenvalue is equal to the Rayleigh quotient; in this case, the bound is proportional to $t[u]^2$, but scaling as $E^2$ for large energy.


\remove{Note that we show that the $\sqrt{E}$ scaling in our theorems
is sharp if the only information
used about the trial function $u$ is the tension $t[u]$ (or $t_s[u]$).
However, numerically, one generally has access to other properties
of $u$ (its normal derivative, whether it is close to $u_\tbox{min}$, etc)
which can give more detailed
statements about the spectrum than merely the distance of $E$ from it.
One example is Still \cite[Thm.~4]{still}
who achieves a bound scaling as $t[u]^2$ via the Kato-Temple inequality.
Another example is the powerful `scaling method' \cite{v+s,que},
which is able to locate many eigenvalues using an operator computed
at a single energy $E$.}

An open problem with practical benefits is the generalization of
these results to Neumann and Robin boundary conditions, and to
multiple subdomains with different trial functions on
each subdomain and least-square errors on artificial boundaries
\cite{descloux,driscoll,timodd} (these are known as Trefftz
or non-polynomial finite element methods).

\section*{Acknowledgments}
The authors are grateful for discussions with Dana Williams, Timo Betcke,
and Chen Hua.
The work of AHB was supported by NSF grant DMS-0811005,
and a Visiting Fellowship to ANU in February 2009 as part of the
{\it ANU 2009 Special Year on Spectral Theory and Operator Theory}.
The work of AH was supported by Australian Research Council Discovery Grant DP0771826. 

\newpage

\appendix 
\section{Proof of Lemma~\ref{l:outgoing}}
\label{a:outgoing}
\bp
For every boundary point $x$, we can find a constant vector field $\mbf{a}_x$ having the property \eqref{outgoing} in a neighbourhood $U_x$ of $x$. (Take a multiple of the vector field $\partial_{z_n}$ in the Euclidean coordinate system used in \eqref{e:Lip}.)
 By compactness we can find a finite  number of such neighbourhoods $U_i = U_{x_i}$, $i = 1, \dots, N$ covering $\Omegab$. We can add to this collection of open sets one additional set $U_0$, whose closure does not meet $\Omegab$, yielding an open cover of $\Omegac$. Let $\phi_i$, $i = 0 \dots N$ be a smooth partition of unity subordinate to this open cover. Then 
$$
\mbf{a} = \sum_{i=1}^N \phi_i \mbf{a}_{x_i}
$$
is a vector field with the required property. 
\ep

\section{Upper bound on distance to spectrum}
\label{a:dist}

\begin{lemma}
Let $\Omega\subset\mathbb{R}^n$ be a bounded domain. Let $E_1$ be the lowest Dirichlet eigenvalue of $\Omega$. Then for any $E > E_1$,
\[
\di{E} \;\le\; C_d E^{1/2}
~,
\qquad C_d = 2 \sqrt{E_1}. 
\]
\label{l:dist}
\end{lemma}
\brmk
The result becomes interesting only for $E> 2(1+\sqrt{2})E_1$. Bounds on
$E_1$ exist as follows.
If $\Omega$ contains a Euclidean ball of radius $r$, then $E_1$ is less than or equal to $E_1(B(0, r))$, which is equal to $j_{n/2 - 1,1}^2/r^2$, where $j_{n/2 - 1,1}$ is the first positive zero of the Bessel function $j_{n/2 - 1}$. For $n=2$ we have $j_{0,1} = 2.4048...$
and for $n=3$ we have $j_{1/2, 1} = \pi$. Also, $E_1$ is greater than or equal to the first eigenvalue of the ball having the same $n$-volume as $\Omega$, by the Faber-Krahn inequality \cite{PolyaSzego}.
\ermk

\bp
Choose a wavevector $\kk \in \Rn$ with $|\kk|^2=E - E_1$, and consider the
trial function $u:\Omega\to\mathbb{C}$ defined by
$u(\xx) := \phi_1(\xx) e^{i\kk\cdot\xx}$, where $\phi_1$ is the
normalized first Dirichlet eigenmode of $\Omega$ with eigenvalue $E_1$.
We calculate,
\[
(\Delta+E)u = 2i \kk\cdot \nabla \phi_1 e^{i\kk\cdot\xx}
~.
\]
Since $u$ is in the domain of $\Omega$, has norm $\|u\|_\lo = 1$, and 
\be
\|(\Delta +E)u\|^2 = 4 \int_\Omega |\kk\cdot\nabla \phi_1|^2 d\xx
\le 4 (E - E_1) E_1 < 4 E E_1,
\label{e:uqm}
\ee
we see that $u$ is an $O(E^{1/2})$ quasimode. 
On the other hand, writing $u=\sum_{j=1}^\infty a_j \phi_j$ we have
\be
\begin{split}
\|(\Delta +E)u\|^2 = \|\sum_j a_j(E-E_j)\phi_j\|^2 = \sum_j |a_j|^2 (E-E_j)^2 \\
\ge \di{E}^2 \sum_j |a_j|^2 = \di{E}^2.
\end{split}
\label{e:usum}
\ee
Combining \eqref{e:uqm} and \eqref{e:usum} completes the proof.
\ep

\newpage

\bibliographystyle{siam} 
\bibliography{alex}

\begin{thebibliography}{10}

\bibitem{Ancona}
{\sc A.~Ancona}, {\em A note on the {R}ellich formula in {L}ipschitz domains},
  Pub. Mathem\`atiques, 42 (1998), pp.~223--237.

\bibitem{backerbdry}
{\sc A.~B{\"a}cker, S.~F{\"u}rstberger, R.~Schubert, and F.~Steiner}, {\em
  Behaviour of boundary functions for quantum billiards}, J. Phys. A, 35
  (2002), pp.~10293--10310.

\bibitem{que}
{\sc A.~H. Barnett}, {\em Asymptotic rate of quantum ergodicity in chaotic
  {E}uclidean billiards}, Comm. Pure Appl. Math., 59 (2006), pp.~1457--88.

\bibitem{incl}
\leavevmode\vrule height 2pt depth -1.6pt width 23pt, {\em Perturbative
  analysis of the {M}ethod of {P}articular {S}olutions for improved inclusion
  of high-lying {D}irichlet eigenvalues}, SIAM J. Numer. Anal., 47 (2009),
  pp.~1952--1970.

\bibitem{mush}
{\sc A.~H. Barnett and T.~Betcke}, {\em Quantum mushroom billiards}, CHAOS, 17
  (2007), p.~043123.
\newblock 13 pages, {\tt nlin.CD/0611059}.

\bibitem{mfs}
\leavevmode\vrule height 2pt depth -1.6pt width 23pt, {\em Stability and
  convergence of the {M}ethod of {F}undamental {S}olutions for {H}elmholtz
  problems on analytic domains}, J. Comput. Phys., 227 (2008), pp.~7003--7026.

\bibitem{timodd}
{\sc T.~Betcke}, {\em A {GSVD} formulation of a domain decomposition method for
  planar eigenvalue problems}, IMA J. Numer. Anal., 27 (2007), pp.~451--478.

\bibitem{gsvd}
\leavevmode\vrule height 2pt depth -1.6pt width 23pt, {\em The generalized
  singular value decomposition and the {M}ethod of {P}articular {S}olutions},
  SIAM J. Sci. Comp., 30 (2008), pp.~1278--1295.

\bibitem{mps}
{\sc Timo Betcke and Lloyd~N. Trefethen}, {\em Reviving the method of
  particular solutions}, SIAM Rev., 47 (2005), pp.~469--491.

\bibitem{chenhua}
{\sc Hua Chen}, {\em Irregular but non-fractal drums, and $n$-dimensional
  {W}eyl conjecture}, Acta Math. Sinica, New Series, 11 (1995), pp.~168--178.

\bibitem{Comech}
{\sc A~Comech}, {\em Cotlar-{S}tein almost orthogonality lemma}, preprint,
  http://www.\-math.tamu.edu/$\sim$comech/papers/CotlarStein/CotlarStein.pdf,
  (2007).

\bibitem{cotlar}
{\sc M~Cotlar}, {\em A combinatorial inequality and its applications to
  $l^2$-spaces}, Rev. Mat. Cuyana, 1 (1955), pp.~41--55.

\bibitem{JK}
{\sc C.~Kenig D.~Jerison}, {\em The inhomogeneous {D}irichlet problem in
  {L}ipschitz domains}, J. Funct. Anal., 130 (1995), pp.~161--219.

\bibitem{descloux}
{\sc J.~Descloux and M.~Tolley}, {\em An accurate algorithm for computing the
  eigenvalues of a polygonal membrane}, Comput. Methods Appl. Mech. Engrg., 39
  (1983), pp.~37--53.

\bibitem{driscoll}
{\sc Tobin~A. Driscoll}, {\em Eigenmodes of isospectral drums}, SIAM Rev., 39
  (1997), pp.~1--17.

\bibitem{fhm}
{\sc L.~Fox, P.~Henrici, and C.~Moler}, {\em Approximations and bounds for
  eigenvalues of elliptic operators}, SIAM J. Numer. Anal., 4 (1967),
  pp.~89--102.

\bibitem{garab}
{\sc P.~R. Garabedian}, {\em Partial differential equations}, John Wiley \&
  Sons Inc., New York, 1964.

\bibitem{gww}
{\sc C.~Gordon, D.~Webb, and S.~Wolpert}, {\em Isospectral plane domains and
  surfaces via {R}iemannian orbifolds}, Invent. Math., 110 (1992), pp.~1--22.

\bibitem{hasselltao}
{\sc Andrew Hassell and Terence Tao}, {\em Upper and lower bounds for normal
  derivatives of {D}irichlet eigenfunctions}, Math. Res. Lett., 9 (2002),
  pp.~289--305.

\bibitem{KSstek}
{\sc J.~R. Kuttler and V.~G. Sigillito}, {\em Inequalities for membrane and
  {S}tekloff eigenvalues}, J. Math. Anal. Appl., 23 (1968), pp.~148--160.

\bibitem{KS}
\leavevmode\vrule height 2pt depth -1.6pt width 23pt, {\em Eigenvalues of the
  {L}aplacian in two dimensions}, SIAM Rev., 26 (1984), pp.~163--193.

\bibitem{litrefftz}
{\sc Z.~C. Li}, {\em The {T}refftz method for the {H}elmholtz equation with
  degeneracy}, Applied Numer. Math., 58 (2008), pp.~131--159.

\bibitem{mclean}
{\sc W.~C.~H. Mc{L}ean}, {\em Strongly elliptic systems and boundary integral
  equations}, Cambridge University Press, 2000.

\bibitem{molerpayne}
{\sc C.~B. Moler and L.~E. Payne}, {\em Bounds for eigenvalues and eigenvectors
  of symmetric operators}, SIAM J. Numer. Anal., 5 (1968), pp.~64--70.

\bibitem{monkwang}
{\sc P.~Monk and D.-Q. Wang}, {\em A least-squares method for the helmholtz
  equations}, Comput. Meth. Appl. Mech. Engrg., 175 (1999), pp.~121--136.

\bibitem{netrusov}
{\sc Yu. Netrusov and Yu. Safarov}, {\em Weyl asymptotic formula for the
  {L}aplacian on domains with rough boundaries}, Commun. Math. Phys., 253
  (2005), pp.~481--509.

\bibitem{PolyaSzego}
{\sc G.~P\'olya and G.~Szego}, {\em Isoperimetric inequalities in mathematical
  physics}, Annals of Mathematics Studies, no. 27, Princeton university press,
  Princeton, NJ, 1951.

\bibitem{rellich}
{\sc Franz Rellich}, {\em Darstellung der {E}igenwerte von {$\Delta u+\lambda
  u=0$} durch ein {R}andintegral}, Math. Z., 46 (1940), pp.~635--636.

\bibitem{SafVas}
{\sc Yu~Safarov and D~Vassiliev}, {\em The Asymptotic Distribution of
  Eigenvalues of Partial Differential Operators}, Translations of Mathematical
  Monographs \#155, American Mathematical Society, Providence, RI, 1996.

\bibitem{saito}
{\sc N.~Saito}, {\em Data analysis and representation on a general domain using
  eigenfunctions of {L}aplacian}, Applied and Computational Harmonic Analysis,
  25 (2008), pp.~68--97.

\bibitem{Seeley}
{\sc R~Seeley}, {\em An estimate near the boundary for the spectral function of
  the {L}aplace operator}, Amer. J. Math., 102 (1980), pp.~869--902.

\bibitem{stein}
{\sc E.~M. Stein}, {\em Harmonic analysis: real-variable methods,
  orthogonality, and oscillatory integrals}, Monographs in Harmonic Analysis,
  Princeton university press, Princeton, NJ, 1993.
\newblock with the assistance of Timothy S. Murphy.

\bibitem{still}
{\sc G.~Still}, {\em Computable bounds for eigenvalues and eigenfunctions of
  elliptic differential operators}, Numer. Math., 54 (1988), pp.~201--223.

\bibitem{tref06}
{\sc Lloyd~N. Trefethen and Timo Betcke}, {\em Computed eigenmodes of planar
  regions}, vol.~412 of Contemp. Math., Amer. Math. Soc., Providence, RI, 2006,
  pp.~297--314.

\bibitem{hakan05}
{\sc H.~E. Tureci, H.~G.~L. Schwefel, P.~Jacquod, and A.~D. Stone}, {\em Modes
  of wave-chaotic dielectric resonators}, Progress in Optics, 47 (2005),
  pp.~75--137.

\bibitem{v+s}
{\sc E.~Vergini and M.~Saraceno}, {\em Calculation by scaling of highly excited
  states of billiards}, Phys. Rev. E, 52 (1995), pp.~2204--2207.

\bibitem{zencyc}
{\sc Steven Zelditch}, {\em Quantum ergodicity and mixing of eigenfunctions},
  in Elsevier Encyclopedia of Mathematical Physics, vol.~1, Academic Press,
  2006, pp.~183--196.
\newblock {\tt arXiv:math-ph/0503026}.

\bibitem{zzw}
{\sc Steven Zelditch and Maciej Zworski}, {\em Ergodicity of eigenfunctions for
  ergodic billiards}, Comm. Math. Phys., 175 (1996), pp.~673--682.

\end{thebibliography}

\end{document}